\documentstyle[12pt]{article}
\title{The divergence of fluctuations for the shape on first passage percolation} 
\author{Yu Zhang\footnote{Research supported by NSF grant DMS-0405150}, Department of Mathematics, University of Colorado}
\date{}
\begin{document}
\baselineskip .20in
\maketitle
\begin{abstract}
Consider the first passage percolation model on ${\bf Z}^d$ for $d\geq 2$.
In this model we assign  independently to each edge the value zero with probability $p$
and the value one with probability $1-p$.
We denote by
$T({\bf 0}, v)$  the passage time from the origin to $v$ for $v\in {\bf R}^d$ and
$$B(t)=\{v\in {\bf R}^d: T({\bf 0}, v)\leq t\}\mbox{ and } G(t)=\{v\in {\bf R}^d: ET({\bf 0}, v)\leq t\}.$$
It is well known that if $p < p_c$, there exists a compact shape $B_d\subset {\bf R}^d$ such that
for all $\epsilon >0$
$$t B_d(1-\epsilon )  \subset {B(t) } \subset tB_d(1+\epsilon )\mbox{ and }
G(t)(1-{\epsilon}) \subset {B(t) } \subset G(t)(1+{\epsilon}) \mbox{ eventually w.p.1.}$$
We denote the fluctuations of $B(t)$ from $tB_d$ and $G(t)$ by
\begin{eqnarray*}
&&F(B(t), tB_d)\!\!=\!\!\inf\left \{l:tB_d\left(1-{l\over t}\right)\subset\!\! B(t)\!\!\subset tB_d\left(1+{l\over t}\right)\right\}\\
&& F(B(t), G(t))\!\!=\!\!\inf\left\{l:G(t)\left(1-{l\over t}\right)\subset B(t)\!\!\subset\!\! G(t)\left(1+{l\over t}\right)\right\}.
\end{eqnarray*}
The means of the fluctuations $E[F(B(t), tB_d]$ and $E[F(B(t), G(t))]$ 
have been conjectured   ranging from divergence to non-divergence for large $d\geq 2$ by physicists. 
In this paper, we show that for all $d\geq 2$ with  a high probability,
the fluctuations $F(B(t), G(t))$ and $F(B(t), tB_d)$
 diverge  with a rate of at least $C \log t$  for some constant $C$.\\

The proof of this argument depends on the linearity between the number of pivotal edges of all minimizing paths
and the paths themselves. This linearity is also independently interesting.\\


\end{abstract}
\section {Introduction of the model and results.}
We consider ${\bf Z}^d, d \geq 2,$ as a graph with edges connecting each pair of vertices
$x=(x_1,\cdots, x_d)$ and $y=(y_1,\cdots, y_d)$ 
with $d(x,y)=1$, where $d(x,y)$ is the distance between $x$ and $y$.
For any two vertex sets $A, B\subset {\bf Z}^d$, the distance between $A$ and $B$ is also
defined by
$$d(A,B)=\min\{d(u,v):u\in A\mbox{ and } v\in B\}.\eqno{(1.0)}$$
We assign independently to each edge the value zero with a probability $p$ or one with probability $1-p$.
More formally, we consider the following probability space. As
sample space, we take $\Omega=\prod_{e\in {\bf Z}^d} \{0,1\},$ points of which
are represented as   {\em configurations}.
Let $P=P_p$ be the corresponding product measure on $\Omega$. The
expectation with respect to $P$ is denoted by $E=E_p$.
For any two vertices $u$ and $v$, 
a path $\gamma$ from $u$ to $v$ is an alternating sequence 
$(v_0, e_1, v_1,...,e_n, v_n)$ of vertices $v_i$ and 
edges $e_i$ in ${\bf Z}^d$ with $v_0=u$ and $ v_n=v$. 
Given a path $\gamma$, we define {\em the passage time} of 
$\gamma$ as 
$$T(\gamma )= \sum_{i=1}^{n} t(e_i).\eqno{(1.1)}$$
For any two sets $A$ and $B$ we define {\em the passage time } from $A$ to $B$
as
$$T(A,B)=\inf \{ T(\gamma): \gamma \mbox{ is a path from some vertex of $A$ to some
vertex in $B$}\},$$
where the infimum takes over all possible finite paths.
A path $\gamma$ from $A$ to $B$ with $t(\gamma)=T(A, B)$ is called the {\em route } of 
$T(A, B)$.
If we focus on a special configuration $\omega$, we may write $T(A, B)(\omega)$
instead of $T(A, B)$. 
When $A=\{u\}$ and $B=\{v\}$ are single vertex sets, $T(u,v)$ is the passage time
from $u$ to $v$. We may extend the passage times over ${\bf R}^d$.
If $x$ and $y$ are in ${\bf R}^d$, we define $T(x, y)=T(x', y')$, where
$x'$ (resp., $y'$) is the nearest neighbor of $x$ (resp., $y$) in
${\bf Z}^d$. Possible indetermination can be dropped by choosing an
order on the vertices of ${\bf Z}^d$ and taking the smallest nearest
neighbor for this order.\\

In particular, the point-point passage time was first introduced by Hammersley and  Welsh (1965)
as follows:
$$a_{m,n}= \inf\{t(\gamma ): \gamma  \mbox{ is a path from $(m,\cdots, 0)$ to  $(n,\cdots 0)$ }\}.
$$
By  Kingman's subadditive argument we have
$$
\lim_{n\rightarrow \infty}{1\over n}  a_{0n} = \mu_p \mbox{ a.s. and in } L_1,\eqno{(1.2)}
$$
and (see Theorem 6.1 in Kesten (1986))
$$\mu_p=0\mbox{ iff } p\geq p_c(d),\eqno{(1.3)}$$
where $p_c=p_c(d)$ is the critical probability for Bernoulli (bond) percolation
on ${\bf Z}^d$ and the non-random constant $\mu_p$ is  called the {\em time constant}. \\

Given a unit vector $x\in {\bf R}^d$, by the same argument in (1.2)
$$
\lim_{n\rightarrow \infty}{1\over n}  T({\bf 0}, nx) = \mu_p(x) \mbox{ a.s. and in } L_1,\eqno{(1.4)}
$$
and 
$$\mu_p(x)=0\mbox{ iff }p\geq p_c.$$

The map $x \rightarrow \mu_p(x)$ induces  a norm on
${\bf R}^d$. The unit radius ball for this norm is denoted by $B_d:=B_d(p)$
and is called {\em the asymptotic shape}. The boundary of $B_d$ is
$$\partial B_d := \{ x \in {\bf R}^d: \mu_p(x)=1\}.$$
 If $p < p_c(d)$, $B_d$ is a compact convex
deterministic set and  $\partial B_d$  is a continuous convex closed
curve (Kesten (1986)).  Define for all $t> 0$
$$B(t):= \{v\in {\bf R}^d, \ T( {\bf 0}, v) \leq t\}.$$
The shape theorem (see Theorem 1.7 of Kesten (1986)) is the well-known result stating that for any
$\epsilon >0$
$$tB_d(1-\epsilon)  \subset {B(t) } \subset tB_d(1+\epsilon )\ \
\mbox{ eventually w.p. 1.}\eqno{(1.5)}$$
In addition to $tB_d$, we can  consider the mean of $B(t)$
$$G(t)=\{v\in {\bf R}^d: ET({\bf 0}, v)\leq t\}.$$
By (1.4), we have
$$tB_d \subset G(t)\eqno{(1.6)}$$
and 
$$G(t)(1-{\epsilon})\subset B(t) \subset G(t)(1+{\epsilon})\mbox{ eventually w.p. 1.}\eqno{(1.7)}$$

The natural or, perhaps the most challenging question in this field (see Kesten (1986) and  Smythe and Wierman (1978)) is to ask ``how fast or how rough"  the boundary is of the set $B(t)$
from the deterministic boundaries $tB_d$ and $G(t)$. This problem has also received a great amount of attention  from 
 statistical physicists.
It is widely conjectured that if $p < p_c(2)$, there exists $\chi(2)=1/3$ 
such that for all $t$ the following probabilities 
$$P\left(tB_2\left (1-{x t^{\chi(2)}\over t}\right) \subset {B(t)} \subset tB_2\left(1+ {xt^{\chi(2)}\over t} \right)\right)\eqno{(1.8)}$$
and 
$$P\left(G(t)\left(1-{x t^{\chi(2)}\over t}\right)  \subset {B(t)} \subset G(t)\left(1+ {xt^{\chi(2)}\over t} \right)\right)\eqno{(1.9)}$$
are close to one or zero for large $x$ or for small $x\geq 0$. In words, the fluctuations of $B(t)$ 
diverge with a rate $t^{1/3}$.
There have been varying in discussions  about the nature of the fluctuations of $B(t)$ for large $d$
ranging from the possible independence of $\chi(d)$ on $d$  
(Kardar and Zhang (1987)) through the picture of $\chi(d)$  decreasing with $d$ 
but always  remaining strictly positive (see Wolf and Kertesz (1987) and Kim and Kosterlitz (1989))
to the possibility that for $d$ above some $d_c$, $\chi(d)=0$ and perhaps the fluctuations do not even diverge (see Natterman and Renz (1988), Haplin-Healy (1989) and Cook and Derrida (1990)).\\

Mathematicians have also made significant efforts in this direction. 
Before we introduce mathematical estimates, let us give a precise definition of the fluctuation of $B(t)$
from some set. For a connected set  $\Gamma$ of ${\bf R}^d$ containing the origin, let
$$\Gamma^+_l= \{v\in {\bf R}^d: d(v, \Gamma)\leq l\} \mbox{ and } \Gamma^-_l=\{v\in \Gamma: d(v,\partial \Gamma)\geq l\}.$$
Note that $\Gamma^-_l\subset \Gamma$ and $\Gamma\subset \Gamma^+_l$. Note also that $\Gamma^-_l$ might be empty
even though $\Gamma$ is non-empty. We say $B(t)$ has a fluctuation 
from $\Gamma$ if
$$F({B(t)},\Gamma)=F(B(t)(\omega), \Gamma)=\inf \{l: \Gamma_l^-\subset B(t)(\omega) \subset \Gamma_l^+\}.\eqno{(1.10)}$$
If we set  $\Gamma =tB_d$ for $d=2$, the conjecture in (1.8) is equivalent to ask if
$$F(B(t), tB_2)\approx  t^{1/3}\mbox{ with a high probability} .$$

When $p\geq p_c(d)$, $B_d$ is unbounded and so is $B(t)$. Also, when $p=0$, there are no fluctuations so we require 
in this paper that
$$0<  P(t(e) =0)=p < p_c (d). \eqno{(1.11)}$$

The mathematical estimates for the upper bound of the fluctuation $F(B(t), \Gamma)$, when $\Gamma=tB_d$ 
and $\Gamma=G(t)$, are quite promising. Kesten (1993) and
Alexander (1993, 1996) showed
that for $p < p_c(d)$ and all $d\geq 2$, there is a constant $C$ such that
$$F(B(t), tB_d)\leq C\sqrt{t} \log t \log t\mbox{ eventually w.p.1}\eqno{(1.12)}$$
and 
$$tB_d \subset G(t) \subset (t+ C_2 \sqrt{t})B_d ,$$
where $\log$ denotes the natural logarithm.
In this paper  $C$ and $C_i$ are always positive  constants that may depend on $p$ or $d$, but not $t$, and their values 
are not significant and  change from 
appearance to appearance.\\

On the other hand, the estimates for the lower bound of the fluctuations are quite unsatisfactory.
Under (1.11) it  seems that the only result for all $d \geq 2$  (see Kesten (1993)) is
$$F(B(t),tB_d)\geq \mbox{ a non-zero constant} \mbox{ eventually w. p. 1}. \eqno{(1.13)}$$
For $d=2$, Piza and Newman (1995) showed that
$$F({B(t)},tB_d)\geq t^{1/8} \mbox{ and } F(B(t), G(t))\geq t^{1/8}\mbox{ eventually w. p. 1}. \eqno{(1.14)}$$

 Clearly,  one of most intriguing  questions in this field is to ask if the fluctuations of $B(t)$
diverge as some statistical physicists believed to be true while others did not. In this paper we answer the conjecture
affirmatively to show that the fluctuations of $B(t)$ always diverge for all $d\geq 2$. We can even tell
that the divergence rate for $B(t)$ is at least $C\log t$. \\

{\bf Theorem 1.} If $0< p < p_c$,  for all $d\geq 3$, $t>0$ and any deterministic set $\Gamma$, 
there exist positive constants $\delta=\delta(p, d)$ $C_1=C_1(p, d)$  such that
$$P(F({B(t)},\Gamma)\geq \delta \log t) \geq 1- C_1 t^{-d+2-2\delta\log p}.  $$

{\bf Remark 1.} If we set $\Gamma=tB_d$ or $\Gamma=G(t)$, together with (1.14),
the fluctuations $F(B(t), tB_d)$ and $F(B(t), G(t))$ 
are at least $\delta \log t$ with a large probability. Also, it follows from this probability estimate that
$$E(F(B(t), \Gamma))\geq C \log t\eqno{(1.15)}$$
for a constant $C=C(p,d)>0$.\\

{\bf Remark 2.} We are unable to  estimate whether $F(tB_d, G(t))$ diverges even though we believe it does.\\

The proof of  Theorem 1 is  constructive. In fact, if $F({B(t)},\Gamma)\leq \delta \log t$,
then we can construct $t^{d-1+2\delta\log p}$ zero paths from $\partial B(t)$ to $\Gamma^+_{\delta \log t}$.
For each such path, we can use the geometric property introduced
in section 2 to show that the path contains  a {\em pivotal} edge defined in section 3.
Therefore, we can construct about $t^{d-1-2\delta\log t}$ pivotal edges. However, in section 3, we can also
show the number
of pivotal edges is of order $t$. Therefore, for a suitable $\delta$ we cannot have as many pivotal edges
as we constructed. The contradiction tells us that $F({B(t)},\Gamma)\geq \delta \log t$.\\

With these estimates for pivotal edges in section 3,
 we can also estimate the number of the total vertices in all routes. This estimate is independently interesting.
For a connected set $\Gamma\subset {\bf R}^d$ with $\alpha_1 tB_d\subset \Gamma \subset \alpha_2 tB_d$
for some constants $0<\alpha_1<1<\alpha_2$, let
$$R_{\Gamma}=\bigcup \gamma_t,\mbox{ where $\gamma_t$ is a route for }T({\bf 0}, \partial \Gamma).$$

{\bf Theorem 2.} If $0< p < p_c(d)$, for all $d \geq 2$ and $t >0$, there exists $C(p, d, \alpha_1, \alpha_2)$ such that
$$E(|R_\Gamma|)\leq C t,$$
where $|A|$ denotes the cardinality of $A$ for some set $A$.\\

{\bf Remark 3.} Kesten (1986) showed that there exists a route for $T({\bf 0}, \Gamma)$
with length $Ct$ for some positive constant $C$. Theorem 2 gives a stronger result with the number of vertices in 
 all routes for $T({\bf 0}, \partial \Gamma)$
 also in order $t$. For quite some time, the author believed that the  routes of $T({\bf 0}, \partial \Gamma)$
resembled a spiderweb
centered at the origin so the number of the vertices in the routes should be of order $t^d$. However,
Theorem 2 negates this assumption.\\ 

{\bf Remark 4.}  Clearly, there might be many routes for the passage time $T({\bf 0}, \partial \Gamma)$.
As a consequence of Theorem 2, each route contains at most $Ct$ vertices. 
Specifically, let
$$M_{x, t}=\sup\{k: \mbox{ there exists a route of $T_\Gamma({\bf 0}, \partial \Gamma)$
containing $k$ edges}\}.$$
If $ 0<p < p_c$,  for all $d\geq 2$ and  $t>0$, 
there exists a positive constant $C=C(p, d)$  such that
$$E(M_{x, t}) \leq Ct.\eqno{(1.16)}$$

{\bf Remark 5.} We may also consider Theorem 2 for a point-point passage time. Let
$$R_{x, t}=\bigcup{ \gamma_n}, \mbox{ where $\gamma_n$ is a route for $T({\bf 0}, xt)$}$$
for a unit vector $x$. We may use the same argument of Theorem 2 to show 
if $0< p< p_c(d)$, for all $d \geq 2$ and $t >0$, there exists a positive constant $C(p, d)$ such that
$$E(|R_{x,t}|)\leq C t.\eqno{(1.17)}$$

{\bf Remark 6.} The condition $p >0$ in Theorem 2 is crucial. As $p\downarrow 0$, the constant $C$ in Theorem 2,
(1.16) and (1.17) may go to infinity. When $p=0$, all edges have to take value one. If we take 
$\Gamma$ as the diamond shape with a diagonal length $2t$ both in vertical and horizontal directions,
it is easy to say that all edges inside the diamond belong to $R_\Gamma$ so
$$|R_\Gamma|=O(t^d).\eqno{(1.18)}$$
This tells us that Theorem 2 will not work when $p=0$.\\

\section{Geometric properties of $B(t)$}

In this section we would like to introduce a few geometric properties for $B(t)$.
We let 
$B'(t)$  be the largest vertex set in $ B(t)\cap {\bf Z}^d$ and $G'(t)$ be the largest vertex set
in  $G(t)\cap {\bf Z}^d$.
Similarly, given a set $\Gamma\subset {\bf R}^d$, we also let 
$\Gamma'$ be the largest vertex set in $\Gamma$.
It is easy to see that
$$ \Gamma'\subset \Gamma\subset \{v+[-1/2, 1/2]^d: v\in \Gamma'\}.  \eqno{(2.0)}$$
As we mentioned in the last section, both $B(t)$ and $G(t)$ are finite as well as $B'(t)$ and $G'(t)$.  
We now  show that $B'(t)$ and $G'(t)$ are  also connected. Here  a set  $A$ is said to be 
connected in ${\bf Z}^d$ if any two vertices of $A$ are connected by 
a path in $A$.\\

{\bf Proposition 1.} $B'(t)$ and $G'(t)$ are  connected.\\

{\bf Proof.} Since $T({\bf 0}, {\bf 0})=0\leq t$, 
${\bf 0}\in B'(t)$. We pick a vertex  $v\in B'(t)$, so $T({\bf 0}, v)\leq t$. This tells us that
there exists a path $\gamma$ such that 
$$T(\gamma) \leq t.$$
Therefore, for any $u\in \gamma$,
$$T({\bf 0}, u)\leq t\mbox{ so } u\in B'(t).$$
This implies that $\gamma \subset B'(t)$, so 
 we know $B'(t)$ is connected. The same argument shows that $G'(t)$ is connected.
 $\Box$\\

Given a finite set $\Gamma$ of ${\bf Z}^d$ we may define  its vertex boundary as follows. For each $v\in \Gamma $,
$v\in \Gamma $ is said to be a boundary vertex of $\Gamma$ if there exists $u\not \in \Gamma$ but $u$
is adjacent to $v$. We denote by $\partial \Gamma$  all boundary vertices of $\Gamma$.
 We also let $\partial_o\Gamma $ be  all vertices not in $\Gamma$, but adjacent to $\partial \Gamma$.
With these definitions, we have the following Proposition.\\

{\bf Proposition 2.} For all $v\in \partial B'(t)$, $T({\bf 0}, v)=t$ and for all $u\in \partial_o B'(t)$
$T({\bf 0}, u)=t+1$.\\

{\bf Proof.} We pick $v\in \partial B'(t)$. By the definition of the boundary, $v\in \partial B'(t)$ so 
$T({\bf 0}, v) \leq t.$ Now we show $T({\bf 0}, v)\geq t$ for all $v\in \partial B'(t)$.
If we suppose that $T({\bf 0}, v)< t$ for some $v\in \partial B'(t)$, then $T({\bf 0}, v) \leq t-1$, since
$T({\bf 0}, v)$ is an integer. Note that $t(e)$ only takes zero and one values
so  there exists $u\in \partial_o B(t)$ and $u$ is adjacent to $v$ such that
$T({\bf 0}, u) \leq t$. This tells us that $u\in B'(t)$. But  we know as we defined that
$$\partial_o B'(t)\cap B(t)=\emptyset.$$
 This contradiction tells us that
$T({\bf 0}, v)\geq t$ for all $v\in \partial B'(t)$.\\

Now we will prove the second part of Proposition 2. We pick a vertex $u\in \partial_o B'(t)$.
Since $u$ is adjacent to $v\in B'(t)$,
$$T({\bf 0}, u)\leq 1+ T({\bf 0}, v)\leq 1+t.$$
On the other hand, any path from ${\bf 0}$ to $u$ has to pass through a vertex of $\partial B'(t)$ before reaching
  $\partial_o B'(t)$.
We denote the vertex by $v$.
As we proved,
$T({\bf 0}, v)= t$.  The passage time of the rest of the path from $v$ to $u$ has to be greater or equal to one,
otherwise, $u\in B'(t)$. Therefore, any path from ${\bf 0}$ to $u$ has a passage time larger or equal to
$t+1$, that is
$$T({\bf 0}, u)\geq t+1.$$
Therefore, $T({\bf 0}, \partial_o B'(t))=t+1.$ $\Box$\\

Given a fixed connected set $\kappa=\kappa_t$ containing the origin, define the 
 event 
$$\{B'(t)=\kappa\}=\{\omega: B'(t)(\omega)=\kappa\}.$$
 
{\bf Proposition 3.} The event that $\{B'(t)=\kappa\}$ only depends on the zeros and ones of the edges of $\kappa \cup \partial_o \kappa $.\\

Proposition 2 for $d=2$ has been proven by Kesten and Zhang (1998). In fact, they gave a precise
structure of $B'(t)$. We may adapt their idea to prove Proposition 2 for $d\geq 3$ by using the plaquette surface
(see the definition in section 12.4 Grimmett (1999)). 
To avoid the complicated definition of the plaquette surface, we would rather give the following direct
proof.\\

{\bf Proof of Proposition 3.} 
Let $\kappa^C$ denote the vertices of ${\bf Z}^d\setminus \kappa$ and 
$$\{\omega(\kappa)\}=\prod_{\mbox{ edge in }\kappa}\{0,1\}\mbox{ and } \{\omega(\kappa^C)\}=\prod_{\mbox{edge in } 
\kappa^C} \{0,1\},$$
where edges in $\kappa$ are the edges whose two vertices  belong to $\kappa\cup \partial_o\kappa$
and the edges in $\kappa^C$ are the other edges.
For each $\omega \in \Omega$, we may rewrite
$\omega$ as 
$$\omega=(\omega(\kappa), \omega(\kappa^C)).$$
Suppose that Proposition 3 is not true, so the zeros and ones in $\omega(\kappa^C)$ can affect the event
$\{B'(t)=\kappa\}$. 
In other words, there exist
 two different $\omega_1,\omega_2 \in \Omega$ with
$$\omega_1=(\omega(\kappa), \omega_1(\kappa^C))\mbox{ and } \omega_2=(\omega(\kappa), \omega_2(\kappa^C))$$
such that 
$$B'(t)(\omega_1)=\kappa\mbox{ but } B'(t)(\omega_2)\neq \kappa.\eqno{(2.4)}$$
From (2.4) there are two cases:\\
(a) there exists $u$ such that $u\in B'(t)(\omega_2)$, but $u\not \in \kappa$.\\
(b) there exists $u$ such that $u\in \kappa$, but $u\not \in B'(t)(\omega_2)$.\\
If (a) holds, 
$$T({\bf 0}, u)(\omega_2)\leq t.\eqno{(2.5)}$$
There exists a path $\gamma$ from ${\bf 0}$ to $u$ such that
$$T(\gamma)(\omega_2)\leq t.\eqno{(2.6)} $$
Since $u\not\in \kappa$, any path from ${\bf 0}$ to $u$ has to pass through $\partial_o \kappa=\partial_o B'(t)(\omega_1)$.
Let $\gamma'$ be the subpath of $\gamma$ from ${\bf 0}$ to $\partial_o \kappa=\partial_o B'(t)(\omega_1)$. Then by
Proposition 2, 
$$T(\gamma')(\omega_1)\geq t+1.$$
Note that the zeros and ones  in both $\omega_2=(\omega(\kappa), \omega_2(\kappa^C)$ and $\omega_1=(\omega(\kappa), \omega_1(\kappa^C)$ are the same
 so
$$t+1\leq T(\gamma')(\omega_1)= T(\gamma')(\omega_2)\leq T(\gamma)(\omega_2).\eqno{(2.7)}$$
By (2.6) and (2.7) (a) cannot hold.\\

Now we assume that (b) holds. Since any path from ${\bf 0}$ to $u$ has to pass through $\partial_oB'(t)(\omega_2)$,
by Proposition 2, 
$$T({\bf 0}, u)(\omega_2)\geq t+1.\eqno{(2.8)}$$
But since $u\in \kappa$ and $B'(t)(\omega_1)=\kappa$, there exists a path $\gamma$ inside $B'(t)(\omega_1)$
from ${\bf 0}$ to $u$ such that
$$T(\gamma)(\omega_1)\leq t.$$
Therefore, 
$$T({\bf 0}, u)(\omega_2) \leq T(\gamma)(\omega_2)\leq t,\eqno{(2.9)}$$
since $\gamma\subset \kappa$ and the  zeros and ones  in both $\omega_2=(\omega(\kappa), \omega_2(\kappa^C)$ and $\omega_1=(\omega(\kappa), \omega_1(\kappa^C)$ are the same.
The contradiction of (2.8) and (2.9) tells us that (b) cannot hold. $\Box$\\

\section{The linearity of a number of pivotal edges}

In this section  we will discuss a fixed value $0<p< p_c$ and a fixed interval open $I_p\subset (0,p_c)$
centered at $p$. First we show that the length of a route from the origin to $\partial B'(t)$
is of order $t$.\\

{\bf Lemma 1.} For a small interval $I_p \subset (0, p_c)$ centered at $p$,
there exist positive constants $\alpha=\alpha(I_p,d)$, $C_1=C_1(I_p, d)$ and $C_2=C_2(I_p, d)$ such that for all $t$ and all $p'\in I_p$,
$$P( \exists \mbox{ a route $\gamma$ from the origin to $\partial B'(t)$ with } |\gamma| \geq \alpha t)
\leq C_1 \exp(-C_2 t).$$

{\bf Proof.} 
By Theorem 5.2 and 5.8  in Kesten (1986)  for all $p' \in I_p$ and for all $t$ there exist $C_3=C_3(I_p,d)$ and $C_4=C_4(I_p,d)$ such that
$$P(2t/3 B_d \not \subset  B(t))\leq C_3 \exp(-C_4 t) \mbox{ and } P(B(t)\not \subset 3t/2 B_d) \leq C_3 \exp(-C_4 t).
\eqno{(3.0)}$$
If we put these two inequalities from (3.0) together, we have for all $p'\in I_p$ and all $t$,
$$P(t/2 B_d)\subset B'(t) \subset 2t B_d) \mbox{ for all large }t)\geq 1-C_3 \exp(-C_4 t).\eqno{(3.1)}$$
On the event of 
$$  \{\mbox{ a route $\gamma$ from the origin to $\partial B(t)$ with } |\gamma| \geq \alpha t\}\cap \{t/2 B_d)\subset B(t) \subset 2t B_d\}$$
we assume that there exists a route $\gamma$ from the origin to some vertex $u\in \partial B'(t)$
with
$$t/2\leq d({\bf 0}, u) \leq 2t\mbox{ and } |\gamma| \geq \alpha t$$
 such that
$$ t=T({\bf 0}, \partial B'(t))= T(\gamma)=T({\bf 0}, u).$$
Therefore, by (3.1)
\begin{eqnarray*}
&&P( \exists \mbox{ a route $\gamma$ from the origin to $\partial B'(t)$ with } T(\gamma)\leq t,|\gamma| \geq \alpha t)\\
&&\leq \sum_{t/2\leq d({\bf 0}, u) \leq 2t} 
P( \exists \mbox{ a route $\gamma$ from the origin to $u$ with } T(\gamma)\leq t,|\gamma| \geq \alpha t)+C_3\exp(-C_4 t).
\end{eqnarray*}
Proposition 5.8 in Kesten (1986) tells us that there exist positive constants $\beta (I_p, d)$, 
$C_5=C_5(I_p, d)$ and $C_6=C_6(I_p, d)$ such that for all $p'\in I_p$ and $t$
\begin{eqnarray*}
&&P(\exists \mbox{ a self avoiding path $\gamma$ from $(0,0)$ to $y$ contains $n$ edges, but
with $T(\gamma) \leq  \beta n$})\\
&&\leq C_5\exp(-C_6 n),\hskip 12cm (3.2)
\end{eqnarray*}
where $n$ is the largest integer less than $t$.
Together with these two observations if we take a suitable $\alpha=\alpha (I_p, d) $, we have for all $p'\in I_p$
$$P( \exists \mbox{ a route $\gamma$ from the origin to $\partial B'(t)$ with } |\gamma| \geq \alpha t)
\leq (2t)^d C_5\exp(-C_6 t).$$
Lemma 1 follows . $\Box$\\

To show Theorems we may concentrate to a ``regular" set satisfying (3.1). Here we give the following precise definition.
Given a deterministic connected finite set $\Gamma=\Gamma_t\subset {\bf R}^d$,
 $\Gamma$ is said to be {\em regular} if there exists $t$ such that
$$1/2 tB_d \subset \Gamma \subset 2 tB_d. \eqno{(3.3)}$$
For a regular set $\Gamma$ we denote by
$$T_\Gamma({\bf 0}, \partial \Gamma)=\inf \{ T(\gamma): \gamma \subset \Gamma'\mbox{ is a path from the origin to some
vertex of $\partial \Gamma'$}\}.$$
Now we try to compute the derivative of $ET_\Gamma({\bf 0}, \partial \Gamma)$ in $p$ for a regular set $\Gamma$.
As a result, we have
$$ET_\Gamma({\bf 0}, \partial \Gamma)=\sum_{i\geq 1} P(T_\Gamma((0,0), \partial \Gamma)\geq i).$$
Note that  $t(e)$ takes zero with probability $p$ and one with  probability $1-p$. 
But we know that the standard Bernoulli random variable takes zero with  probability $1-p$ and one with  probability $p$.
Hence we have to define an increasing event (see the definition in section 2 of
Grimmett (1999)) in reverse.
An event ${\cal A}$ is said to be increasing if
$$1-I_{\cal A}(\omega) \leq 1-I_{\cal A}(\omega') \mbox{ whenever } \omega\leq \omega',$$
where $I_{\cal A}$ is the indicator of ${\cal A}$.
 Note that $\Gamma$ is a finite set, so ${d E T_\Gamma({\bf 0}, \partial \Gamma)\over dp}$ exists.
We have
$${d E T_\Gamma({\bf 0}, \partial \Gamma)\over dp}  =\sum_{i\geq 1} {d P(T_\Gamma({\bf 0}, \partial \Gamma)\geq i)\over dp}.\eqno{(3.5)}$$
 Note that 
$$\{ T_\Gamma({\bf 0}, \partial \Gamma)\geq i\}$$
is decreasing  so
by Russo's formula
$${d E_p T_\Gamma({\bf 0}, \partial \Gamma)\over dp}  =-\sum_{i\geq 1}\sum_{e\in \Gamma} P(\{T_\Gamma({\bf 0}, \partial \Gamma)\geq i\}(e)),\eqno{(3.6)}$$
where $\{T_\Gamma({\bf 0}, \partial \Gamma)\geq i\}(e)$ is the event that $e$ is a pivotal for $\{T_\Gamma({\bf 0}, \partial \Gamma)\geq i\}$.
In fact, given  a configuration $\omega$,  $e$ is said to be a pivotal
edge for $\{T_\Gamma({\bf 0}, \partial \Gamma)(\omega)\geq i\}$ if $t(e)(\omega)=1$ and
$$T_\Gamma({\bf 0}, \partial \Gamma)(\omega')=i-1$$
where $w'$ is the configuration that $t(b)(\omega)=t(b)(\omega')$ for all edges $b\in \Gamma$ except $e$ and
$t(e)(\omega')=0$. The event  $\{T_\Gamma({\bf 0}, \partial \Gamma)\geq i\}(e)$ is equivalent to the event that
there exists a route of $T_\Gamma({\bf 0}, \partial \Gamma)$ with passage time $i$
passing  through $e$ and $t(e)=1$.
With this observation,
\begin{eqnarray*}
&&{d E T_\Gamma({\bf 0}, \partial \Gamma)\over dp}  \\
&=&-\sum_{i\geq 1}\sum_{e\in \Gamma} P(\exists \mbox{ a route
of $T_\Gamma({\bf 0}, \partial \Gamma)$ passing through $e$ with $T_\Gamma({\bf 0}, \partial \Gamma)=i$ and $t(e)=1$})\\
&&=-\sum_{e\in \Gamma} P(\exists \mbox{ a route
of $T_\Gamma({\bf 0}, \partial \Gamma)$ passing through $e$ and $t(e)=1$}).
\end{eqnarray*}
Let $K_\Gamma$ be the number of edges $\{e\}\subset \Gamma'$ such that a route from the origin to $\partial \Gamma'$
passes through $e$ and $t(e)=1$.
We have
$$ -{d E T_\Gamma({\bf 0}, \partial \Gamma)\over dp} =E(K_\Gamma).\eqno{(3.7)}$$

Now we give an upper bound for $E(K_\Gamma)$ by giving an upper bound for 
$-{d E T_\Gamma({\bf 0}, \partial \Gamma)\over dp}$.  
Before doing that, we shall define the route length for 
$T_\Gamma({\bf 0}, \partial \Gamma)$
by
$$N_\Gamma(\omega)=\min\{k: \mbox{ there exists a route of $T_\Gamma({\bf 0}, \partial \Gamma)(\omega)$
containing $k$ edges}\}.$$
We  show that the size of $N_\Gamma$ cannot be more than $Ct$ for some constant $C$.\\

{\bf Lemma 2.} For a regular set $\Gamma$ and the interval $I_p$,
there exist positive constants $C_i=C_i(I_p,d)$ $(i=1,2,3$) such that for all $p'\in I_p$ and $t$
$$P( N_\Gamma\geq C_1 t)\leq C_2\exp(-C_3t).$$

{\bf Proof.} We follow the proof of Theorem 8.2 in  Smythe and Wierman (1979).
Let $\omega+r$ denote the time state of the lattice obtained by
adding the $r$ to $t(e)$ for each edge $e$. It follows from the definitions of the passage time and $N_\Gamma$
$$T_\Gamma({\bf 0}, \partial \Gamma)(\omega+r) \leq T_\Gamma({\bf 0}, \partial \Gamma)(\omega)+r N_\Gamma(\omega).\eqno{(3.8)}$$
If we take a negative $r$ in (3.8), we have
$$N_\Gamma (\omega)\leq {T_\Gamma({\bf 0}, \partial \Gamma )(\omega+r)-T_\Gamma({\bf 0}, \partial \Gamma)\over r}.\eqno{(3.9)}$$
Note that $\Gamma$ is regular, so 
$\Gamma\subset 2t B_d$. If we denote by $L$ the segment from the origin to $\partial (2t B_d)$ along the $X$-axis,
then $L$ has to go through $\partial \Gamma$ somewhere since $\Gamma \subset 2t B_d$. Therefore,
$${-T_\Gamma ({\bf 0}, \Gamma)\over r}\leq -{T(L)\over r} \leq -{2 t\over \mu r}.\eqno{(3.10)}$$
If we can show that for some $r < 0$, there exist constants  $C_4=C_4(I_p,d)$  and $C_5=C_5(I_p,d)$ 
such that for all $p'\in I_p$  and all $t$
$$P(T_\Gamma({\bf 0}, \partial \Gamma)(\omega+r)\leq 0) \leq C_4\exp(-C_5t),\eqno{(3.11)}$$
then by (3.9) and (3.10), Lemma 2 holds.
Therefore, to show Lemma 2, it remains to show (3.11).
Note that $\Gamma$ is a finite connected set so for each $\omega$ there exists $x=x(\omega)\in \partial \Gamma$ such that
$$T_\Gamma({\bf 0}, \partial \Gamma)(\omega+r)=T_\Gamma({\bf 0}, x)(\omega+r)\geq T({\bf 0}, x)(\omega+r).$$
Since $x\in \partial \Gamma$ and $\Gamma$ is regular, then
$$t/2\leq d({\bf 0},x)\leq 2t.$$
We have 
$$P(T_\Gamma({\bf 0}, \partial \Gamma)(\omega+r) \leq 0)
\leq P(T({\bf 0}, x(\omega))(\omega+r) \leq 0)
\leq \sum_{t/2\leq d({\bf 0},y)\leq 2t} P(T({\bf 0}, y)(\omega+r)\leq 0).\eqno{(3.12)}$$
Therefore, by (3.2) and (3.12) we take $\beta$ and $|r|$ small with  $r< 0$  and $\beta > |r| >0$ to obtain for all $p'\in I_p$
\begin{eqnarray*}
&& \sum_{t/2\leq d({\bf 0}, y)\leq 2t}P(T({\bf 0}, y)(\omega+r)\leq 0) \\
&&\leq \sum_{t/2\leq d({\bf 0}, y)\leq 2t}P(\exists \mbox{ a self avoiding path $\gamma$ from ${\bf 0}$ to $y$ contains $n$ edges,}\\
&& \hskip 3cm \mbox{ but  with $T(\gamma)(\omega) \leq 2\beta n$})\\
&&\leq C_6t^{d}\exp(-C_7 t),\hskip 11cm (3.13)
\end{eqnarray*}
where $n$ is the largest integer less than $t$.
Therefore, (3.11) follows from (3.13). $\Box$\\

With Lemma 2 we are ready to give  an upper bound for   $-{d E T_\Gamma({\bf 0}, \partial \Gamma)\over dp}$.\\

{\bf Lemma 3.} For a regular set $\Gamma$ there exists a constant $C(I_p, d)$ such that for all $p'\in I_p$ and $t$
$$-{d E T_\Gamma({\bf 0}), \partial \Gamma)\over dp} \leq Ct.$$

{\bf Proof.} 
We assign  $s(e)\geq t(e)$ either zero or one independently from
edge to edge with probabilities  $p-h$ or $1-(p-h)$ for a small number $h>0$, respectively. 
With this definition,
\begin{eqnarray*}
&&P(s(e)=1, t(e)=0)=P(s(e)=1)-P(s(e)=1, t(e)=1)\\
&&=P(s(e)=1)-P(t(e)=1)=1-(p-h)-(1-p)=h.\hskip 4cm {(3.14)}
\end{eqnarray*}
Let $\gamma^t$ be a route for $T_\Gamma^t({\bf 0}, \Gamma)$ with time state $t(e)$ and let
$\gamma^s$ be a route for $T_\Gamma^s({\bf 0}, \Gamma)$  with time state $s(e)$.
Here we pick $\gamma^t$ such that 
$$|\gamma^t|= N_\Gamma.$$
For each edge $e\in \gamma^t$, if $t(e)=1$, then $s(e)=1$. If $t(e)=0$ but $s(e)=1$, we just add one for
this edge. Therefore,
$$T_\Gamma^s({\bf 0}, \Gamma)\leq T(\gamma^t) + \sum_{e\in \gamma^t} I_{(t(e)=0, s(e)=1)}.\eqno{(3.15)}$$
Clearly, $\gamma^t$ may not be unique, so we select a route from these  $\gamma^t$ in a unique way.
We still write $\gamma^t$ for the unique route without loss of generality. 
By (3.15) and this selection
$$ET_\Gamma^s({\bf 0}, \Gamma)\leq ET(\gamma^t) + \sum_{\beta}\sum_{e\in \beta } P(t(e)=0, s(e)=1,\gamma^t=\beta),\eqno{(3.16)}$$
where the first sum in (3.16) takes over all possible paths $\beta$ from ${\bf 0}$ to $\partial \Gamma'$.
Let us estimate 
$$\sum_{\beta}\sum_{e\in \beta } P(t(e)=0, s(e)=1,\gamma^t=\beta).$$
Since $\Gamma$ is regular, the longest path from $(0,0)$ to $\partial \Gamma'$ is less than $(2t)^d$. By Lemma 2,
there exist $C_1=C_1(I_p,d)$, $C_2(I_p, d)$ and $C_3(I_p, d)$ such that
\begin{eqnarray*}
&&\sum_{\beta}\sum_{e\in \beta } P(t(e)=0, s(e)=1,\gamma^t=\beta)\\
&\leq &\sum_{|\beta|\leq C_1 t}\sum_{e\in \beta } P(t(e)=0, s(e)=1,\gamma^t=\beta)+  C_2 t^d \exp(-C_3t).\hskip 3cm (3.17)
\end{eqnarray*}
Note that the value of $s(e)$ may depend on the value of $t(e)$, but not the other values of $t(b)$ for $b\neq e$, so
by (3.14)
\begin{eqnarray*}
&&P(t(e)=0, s(e)=1,\gamma^t=\beta)=P(s(e)=1\,\, |\,\, t(e)=0, \gamma^t=\beta)P(t(e)=0,\gamma^t=\beta)\\
&&\leq 
P(s(e)=1\,\, |\,\, t(e)=0)P(\gamma^t=\beta)=hp^{-1}P(\gamma^t=\beta).\hskip 4cm (3.18)
\end{eqnarray*}
By  (3.18) we have
$$\sum_{|\beta|\leq C_1}\sum_{e\in \beta } P(t(e)=0, s(e)=1,\gamma^t=\beta)\leq \sum_{|\beta|\leq C_1t}\sum_{e\in \beta } hp^{-1}P(\gamma^t=\beta)\leq C_1htp^{-1}.\eqno{(3.19)}$$

By  (3.17) and (3.19) there exists $C_4=C_4(I_p, d)$ such that
$$E(T_\Gamma^s({\bf 0}, \partial \Gamma))\leq E(T_\Gamma^t({\bf 0}, \partial \Gamma)+C_4th.\eqno{(3.20)}$$
If we set
$$f(p)= E(T_\Gamma({\bf 0}, \partial \Gamma)) \mbox{ for time state $t(e)$ with $P(t(e)=0)=p$},$$
then  by (3.20)
$$-{df(p)\over dp} =\lim_{h\rightarrow 0} -{f(p-h)-f(p)\over -h}\leq C_4t.\eqno{(3.21)}$$ 
Therefore, we have
$$-{d E T_\Gamma({\bf 0}, \partial \Gamma)\over dp}=-{df(p)\over dp}\leq C_4t.\eqno{(3.22)}.$$
Therefore, Lemma 3 follows from  (3.22). $\Box$\\

Together  with (3.7) and Lemmas 3, we have the following proposition.\\

{\bf Proposition  4.} If $0<p < p_c$, then for a regular set $\Gamma$
there exists a constant $C=C(p)$ such that
$$ EK_\Gamma\leq Ct.$$

\section{ Proof of  Theorem 1}
In this section, we only show Theorem 1 for $d=3$. The same proof for $d>3$ can be adapted directly.
Given a fixed  set $\Gamma\subset {\bf R}^3$ defined in section 2, $\Gamma'\subset {\bf Z^3}$ is
 the largest vertex set in $\Gamma$, where 
$$\Gamma'\subset \Gamma\subset \{v+[1/2,1/2]^3: v\in \Gamma'\}.\eqno{(4.0)}$$
Suppose that there exists a deterministic set $\Gamma$ such that
$$F(B(t), \Gamma) \leq \delta \log t .\eqno{(4.1)} $$
(4.1) means that 
$$\Gamma^-_{\delta \log t} \subset B(t) \subset \Gamma^+_{\delta \log t}, $$
where
$$\Gamma^+_l= \{v\in {\bf R}^3: d(v, \Gamma)\leq l\} \mbox{ and } \Gamma^-_l=\{v\in \Gamma: d(v,\partial \Gamma)\geq l\}.\eqno{}$$

We first show that if $\Gamma^+_{\delta \log t}$ does not satisfy the regularity condition in (3.3), then the probability
of the event in (4.1) is exponentially small.
We assume that 
$$\Gamma^+_{\delta \log t} \not \subset 2t B_d.\eqno{(4.2)}$$
If 
$$F(B(t), \Gamma) \leq  \delta \log t\mbox{ with } \delta \log t < {t\over 3}, \eqno{(4.3)}$$
then we claim that 
$$B(t) \not \subset {3t \over 2} B_d. \eqno{(4.4)}$$
To see (4.4), note that
$$B(t)\subset {3t\over 2} B_d \mbox{  implies that } \Gamma^+_{\delta \log t}\subset {2t} B_d.\eqno{(4.5)}$$
Therefore, (4.4) follows from (4.2).
Under (4.2),  by (3.0) there exist $C_1(p,d)$ and $C_2(p, d)$ such that
$$P(F(B(t), \Gamma) \leq  \delta \log t) \leq C_1 \exp(-C_2 t).\eqno{(4.6)}$$
Similarly, if we assume that 
$(t/2) B_d \not\subset \Gamma^+_{\delta \log t} $ for a set $\Gamma$, we have
$$P(F(B(t), \Gamma) \leq  \delta \log t) \leq C_1 \exp(-C_2 t).\eqno{(4.7)}$$
With (4.6) and  (4.7), if $\Gamma^+_{\delta\log t}$ does not satisfy the regularity condition in (3.3),
$$P((F(B(t), \Gamma) \leq  \delta \log t) \leq C_1 \exp(-C_2 t).\eqno{(4.8)}$$

Now we focus on $\Gamma^+_{\delta\log t}$ satisfying (3.3). 
We need to show that under (4.1) there are of order  $t^2$ disjoint  
zero paths from
$\partial B(t)$ to $\Gamma^+_{\delta \log t}$. To accomplish this,
let $S_{mt}$ denote a sphere with the center at the origin and a radius $tm$ for small but positive number $m$.
Then by (3.1) for a suitable $m>0$ 
$$P( S_{mt} \subset B(t) \subset 2t B_d)\geq 1-C_1\exp(-C_2 t).\eqno{(4.9)}$$
Here we select the sphere $S_{mt}$ without a special purpose since the sphere is easy to describe.
For each $s\in \partial S_{mt}$, let $L_s$ be the normal line passing though $s$, that is 
the line  orthogonal to the tangent plane of $S_{mt}$ at $s$.
We denote the cylinder with center at $L_s$ by (see Fig.1)
\begin{figure}\label{F:graphG}
\begin{center}
\setlength{\unitlength}{0.0125in}%
\begin{picture}(200,200)(67,800)
\thicklines
\put(150,900){\oval(208,208)}
\put(150,900){\circle{1000}{$S_{mt}$}}
\put(160,920){\line(1,2){50}}
\put(140,920){\line(1,2){50}}
\put(173,905){\line(1,2){55}}
\put(187,975){\line(0,1){10}}
\put(187.5,975){\line(0,1){10}}
\put(188,975){\line(0,1){10}}
\put(188.5,975){\line(0,1){10}}
\put(189,975){\line(0,1){10}}
\put(187,985){\line(1,0){10}}
\put(197,985){\line(0,1){10}}
\put(197,995){\line(1,0){10}}
\put(207,995){\line(0,1){10}}
\put(180,988){$\small{v_i}$}
\put(172,980){$\small{e_{v_i}}$}
\put(198,988){$\gamma_{v_i}$}
\put(210,1025){$L_{s_i}$}
\put(152,912){${s_i}$}
\put(152,940){{$T_{s_i}$}}
\put(187,975){\line(1,0){20}}
\put(207,975){\line(0,-1){20}}
\put(207,955){\line(1,0){10}}
\put(217,955){\line(0,-1){30}}
\put(217,925){\line(1,0){10}}
\put(227,925){\line(0,-1){60}}
\put(227,865){\line(-1,0){10}}
\put(217,865){\line(0,-1){40}}
\put(217,825){\line(-1,0){50}}
\put(167,825){\line(0,-1){10}}
\put(167,815){\line(-1,0){40}}
\put(127,815){\line(0,1){10}}
\put(127,825){\line(-1,0){40}}
\put(87,825){\line(0,1){30}}
\put(87,855){\line(-1,0){20}}
\put(67,855){\line(0,1){40}}
\put(67,895){\line(1,0){10}}
\put(77,895){\line(0,1){50}}
\put(77,945){\line(1,0){10}}
\put(87,945){\line(0,1){40}}
\put(87,985){\line(1,0){60}}
\put(147,985){\line(0,-1){10}}
\put(147,975){\line(1,0){40}}
\put(260,900){${\partial (\Gamma^+_{\delta \log t})'}$}
\put(190,900){${\partial B'(t)}$}

\end{picture}
\end{center}
\caption{The graphs: $S_{mt}$, $\partial B'(t)$, $\Gamma_{\delta\log t}^+$,
 the cylinder $T_{s_i}$ with the center at $L_{s_i}$, pivotal edge $e_{v_i}$,  zigzag path $\gamma_{v_i}$ from $v_i$ to $\partial (\Gamma_{\delta\log t}^+)'$.}
\end{figure}
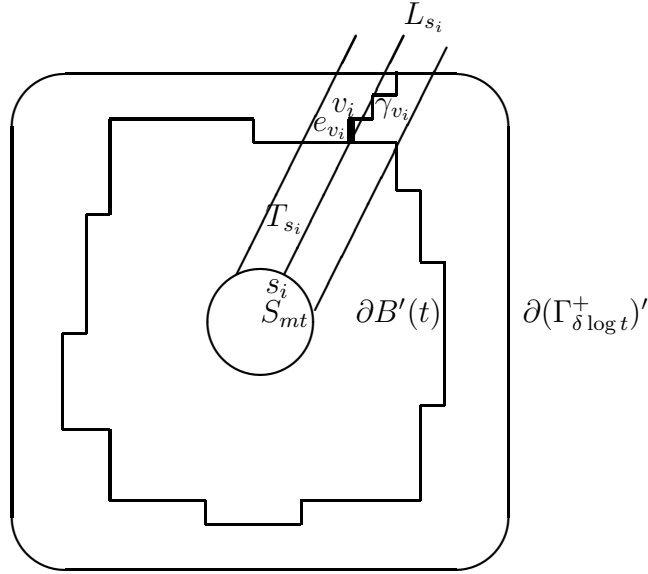 

$$T_s(M) =\{(x,y,z)\in {\bf R}^3: d((x,y,z), L_s)\leq M\}\mbox{ for some constant }M>0.$$
Now we work on the regular polyhedron with $ct^2$ faces embedded on $S_{mt}$, where $ct^2$ is an integer and
 $c=c(m, M)$ is a small number such that
the radius of each face of the regular polyhedron is larger than $M$. 
We denote the center of each face in the regular polyhedron
by $\{s_i\}_{i=1}^{ct^2}$.
By this construction, we have at least $ct^2$ disjoint cylinders $\{T_s(M)\}$.
We denote them  by $\{T_{s_i}(M)\}_{i=1}^{ct^2}$.\\

For each $s_i$, we may take $M$ large such that there exists a path $\gamma_{s_i}\subset Z^d\cap T_{s_i}(M)$
from some vertex of ${\bf Z}^3$ in $S_{mt}$ to $\infty$.
To see the existence of such a path if $L_{s_i}$ is the ray going along the coordinate
axis, we just use $L_{s_i}$ as the path. If it is not, we can construct a zigzag path in 
${\bf Z}^d$ next to $L_{s_i}$ from $s_i$ to $\infty$ (see Figure 1).
In fact, we may take $M =2$ to keep our zigzag path  inside $T_{s_i}(M)$.
For simplicity, we use $T_{s_i}$ to denote $T_{s_i}(2)$.
There might be many such zigzag paths, so we just select one in a unique manner.
We denote by $u_i\in {\bf Z}^d$ with $d(u_i, s_i)\leq 2$ the initial vertex in $\gamma_{s_i}$.
Since $\gamma_{s_i}$ is next to $L_{s_i}$,
for any  point $x$ on the
ray $L_{s_i}$, there is $v$ in $\gamma_{s_i}$ with 
$d(v,y)\leq 2$. Furthermore, by a simple induction we conclude that
$$\mbox{the number of vertices from $u_i$ to $v$ along $\gamma_{s_i}$ is less than $2d(s_i,x)$}.\eqno{(4.10)}$$
(4.10) tells us that the length of $\gamma_{s_i}$  is linear to the length of $L_{s_i}$.
 Since $\gamma_{s_i}$ is from $S_{mt}$ to $\partial (\Gamma^+_{\delta\log t})'$, it has to come
outside of $B'(t)$ from its inside. 
Let $v_i$ be the vertex
in $\partial_o B'(t)$ and $\gamma _{v_i}$ be the piece of $\gamma_{s_i}$ outside $B'(t)$
from $v_i$ to $\partial (\Gamma^+_{\delta \log t})'$ (see Figure 1).
On $F(B(t),\Gamma)\leq \delta \log t$ for a regular $\Gamma$, we know 
$$B'(t)\subset (\Gamma_{\delta\log t}^+)'.$$
Therefore, by our construction (see Figure 1)
$$\gamma_{v_i} \subset (\Gamma_{\delta\log t}^+)'.\eqno{(4.11)}$$
Also, by our special construction in (4.10), we have
$$|\gamma_{v_i}|\leq 2 \delta \log t.\eqno{(4.12)}$$

When $B'(t)=\kappa$ for a fixed vertex set $\kappa$, 
then $\gamma_{v_i}$ is a  fixed path from $\partial_o \kappa$ to $\partial (\Gamma^+_{\delta \log t})'$
with a length less than $2 \delta \log t$.
Therefore, on $B'(t)=\kappa$ 
$$P( \gamma_{v_i} \mbox{ is  a zero path})\geq p^{2\delta \log t}.\eqno{(4.13)}$$
We say $T_{s_i}$ is good if there exists such a zero path $\gamma_{v_i}$.
On $B'(t)=\kappa$, 
let $M(\Gamma, \kappa)$ be the number of such good cylinders $T_{s_i}$. By (4.13), we have
$$EM(\Gamma, \kappa)\geq (ct^2)p^{2\delta \log t}=ct^{2-2\delta\log p}.\eqno{(4.14)}$$
On $B'(t)=\kappa$,  note that the event that $T_{s_i}$ is good depends on zeros and ones of
the edges inside $T_{s_i}$, but outside of $\partial_o\kappa$. Note also that $T_{s_i}$ and $T_{s_j}$ are disjoint
for $i\neq j$, so by a standard Hoeffding  inequality there exist $C_1=C_1(p,d)$ and $C_2=C_2(p,d)$ such that
$$P(M(\Gamma, \kappa) \leq ct^{(2-2\delta\log p)}/2) \leq C_1 \exp(-C_2 ct^{-2+2\delta\log p}).\eqno{(4.15)}$$
We denote by 
$${\cal D}(\Gamma, \kappa)=\{ M(\Gamma, \kappa) \geq ct^{(2+2\delta\log p)}/2\}.$$
Note that 
$${\cal D}(\Gamma, \kappa)\mbox{ only depends zeros and ones outside }\partial_o\kappa. \eqno{(4.16)}$$
By Proposition 3 and  (4.16),
$$\{B'(t)=\kappa\} \mbox{ and }{\cal D}(\Gamma, \kappa)\mbox{ are independent. }  \eqno{(4.17)}$$
By Proposition 2, any route from $(0,0,0)$ to $v_i$ in  $B'(t)\cup \partial_o B'(t)$ has a passage time $t+1$.
We just pick one from these routes and denote it by $\gamma({\bf 0}, v_i)$.
On $F(B(t),\Gamma) \leq \delta \log t$ if $T_{s_i}$ is good, there exists a zero path $\gamma_{v_i}$ from $v_i$ to $\partial (\Gamma^+_{\delta \log t})'$.
This implies that there exists a path 
$$\gamma({\bf 0}, \partial \Gamma^+_{\delta \log t})=\gamma({\bf 0}, v_i)\cup \gamma_{v_i}$$
from $(0,0,0)$ to $\partial (\Gamma^+_{\delta \log t})'$
with a passage time $t+1$ and the path passes through the edge adjacent $v_i$ between $\partial B(t)$ and $\partial_o B(t)$.
On the other hand,  note that any path from the origin to $\partial (\Gamma^+_{\delta \log t})'$ 
has to pass through $\partial_o B'(t)$ first, so by Proposition 2 it has to spend at least passage time $t+1$. 
Therefore, if we denote by $e_{v_i}$ the edge adjacent $v_i$ from $\partial B(t)$ to
$\partial_o B(t)$, then the path $\gamma({\bf 0}, \partial \Gamma^+_{\delta \log t})$  with  passage time $T((0,0,0), \partial \Gamma^+_{\delta\log t})$ passes through $e_{v_i}$
and $t(e_{v_i})=1$. 
By (4.11) and
$$B'(t)\subset (\Gamma_{\delta\log t}^+)',$$
the path $\gamma({\bf 0}, \partial \Gamma^+_{\delta \log t})$ has to stay inside $(\Gamma^+_{\delta \log t})'$.
These observations tell us that $e_{v_i}$ is a pivotal edge for  $T_{\Gamma^+_{\delta \log t}}((0,0), \Gamma^+_{\delta \log t})$.
Therefore, on $F(B(t),\Gamma) \leq \delta \log t$ if $T_{s_i}$ is good,
$$T_{s_i}\mbox{ contains at least one pivotal edge for $T_{\Gamma^+_{\delta \log t}}((0,0), \Gamma^+_{\delta \log t}$)}.\eqno{ (4.18)}$$

With these preparations we now show Theorem 1. \\

{\bf Proof of Theorem 1.}
If $\Gamma^+_{\delta \log t}$ is not regular, 
$$1/2tB_d \not\subset \Gamma^+_{\delta \log t} \mbox{ or } \Gamma^+_{\delta \log t} \not\subset 2t B_d,$$ 
by (4.8) there are $C_1=C_(p,d)$ and $C_2(p, d)$ such that
$$P_p(F(B(t), \Gamma)\leq \delta \log t)\leq C_1\exp(-C_2 t).\eqno{(4.19)}$$
Now we only need to focus on a regular  $\Gamma_{\delta \log t}^+$.
$$P_p(F(B(t), \Gamma)\leq \delta \log t)
= \sum_{\Gamma}P(F(B(t), \Gamma )\leq \delta \log t,B'(t)=\kappa),\eqno{(4.20)}
$$
where the sum takes over all possible sets $\kappa$.
For each fixed $\kappa$, by (4.17) and (4.15) there exist $C_3=C_3(p,d)$ and $C_4=C_4(p,d)$ such that
\begin{eqnarray*}
&&\sum_{\kappa}P(F(B(t), \Gamma )\leq \delta \log t,B'(t)=\kappa)\\
\leq&& \sum_{\kappa}P(F(B(t), \Gamma )\leq \delta \log t,B'(t)=\kappa, {\cal D}(\Gamma, \kappa))
+C_3\exp(-C_4 t^{2+2\delta \log p}).
\end{eqnarray*}
By (4.18)
$$\!\!\!\!\!\!\!\!\!\!\!\! \sum_{\kappa}P\left(F(B(t), \Gamma )\leq \delta \log t,B'(t)=\kappa, {\cal D}(\Gamma, \kappa)\right)
\leq \sum_{\kappa}P\left(K_{\Gamma ^+_{\delta \log t}}\geq { ct^{2+2\delta\log p}\over 2} ,B'(t)=\kappa \right).\eqno{(4.21)}$$
We combine (4.20)-(4.21) together to have 
$$P(F(B(t), \Gamma)\leq \delta \log t)\leq P\left(K_{\Gamma^+_{\delta \log t}}\geq {ct^{2+2\delta\log p}\over 2}\right)+C_5\exp(-C_6 t^{2+2\delta \log p})\eqno{(4.22)}$$
for $C_5=C_5(p, d)$ and $C_6=C_6(p,d)$.
By  Markov's inequality and Proposition 4, if we select a suitable $\delta>0$, for a regular $\Gamma$ there exists
$C_7=C_7(p, d)$ such that
$$P(F(B(t), \Gamma)\leq \delta \log t)
\leq C_7t^{-1-2\delta\log p}.\eqno{(4.23)}$$
Theorem  1 follows from   (4.19) and (4.23).

\section { Proof of Theorem 2.} 
Since $\Gamma$ is regular, by Proposition 4,
$$E K_\Gamma \leq Ct.\eqno{(5.1)}$$
By (5.1) for a large positive number $M$,
$$E|R_\Gamma|\leq E(|R_\Gamma|; |R_\Gamma |\geq M K_\Gamma)+MCt\eqno{(5.2)}$$
Now we estimate $E(|R_\Gamma|; |R_\Gamma |\geq M K_\Gamma)$
by using the method of renormalization in Kesten and Zhang (1990).
We define, for integer $k\geq 1$ and $u\in {\bf Z}^d$, the cube
$$B_k(u)=\prod_{i=1}^d[ku_i, ku_i+k)$$
with lower left hand corner at $ku$ and fattened $R_\Gamma$ by
$$\hat{R}_\Gamma(k)=\{u\in {\bf Z}^d: B_k(u)\cap R_\Gamma \neq \emptyset \}.$$
By our definition, 
$$ |\hat{R}_\Gamma(k)|\geq {|R_\Gamma|\over k^d}.\eqno{(5.3)}$$

For each cube $B_k(u)$, it has at most $4^d$ neighbor cubes, where we count its diagonal neighbor cubes.
We say these neighbors are connected to $B_k(u)$ and denote by $\bar{B}_k(u)$ the vertex set of $B_k(u)$ and
all of its neighbor cubes.
Note that, by the definition, $R_\Gamma$ is a connected set that contains the origin, so $\hat{R}_\Gamma(k)$ is also
connected in the sense of the connection of two of its diagonal vertices.
If $\Gamma$ is regular, then $|R_\Gamma|\geq t/2$. By (5.3), we have
$$P(|R_\Gamma|\geq Mt)=\sum_{m\geq Mt/(2k^d)} P(|R_\Gamma|\geq MK_\Gamma, |\hat{R}_\Gamma(k)|=m).\eqno{(5.4)}$$
We say a cube $B_k(u)$ for $u\in  \hat{R}_\Gamma(k)$ is {\em bad}, if there does not exist an edge $e\in \bar{B}_k(u)\cap R_\Gamma$ such that $t(e)=1$.
Otherwise, we say the cube is {\em good}.
Let ${\cal B}_k(u)$ be the event that  $B_k(u)$ is bad 
and let  $D_\Gamma$ be the number of bad cubes $B_k(u)$ for  $u\in \hat{R}_\Gamma$.\\

If ${\cal B}_k(u)$ occurs, there is a zero path $\gamma_k\subset R_\Gamma $ from $\partial B_k(u)$ to $\partial \bar{B}_k(u)$. By Theorem 5.4 in Grimmett (1999), there exist $C_1=C_1(p,d)$ and $C_2=C_2(p,d)$ 
such that  for fixed $B_k(u)$
$$P({\cal B}_k(u))\leq C_1\exp(-C_2k).\eqno{(5.5)}$$
On $\{|\hat{R}_\Gamma (k)|=m, |R_\Gamma|\geq M K_\Gamma\}$, if $2(4k)^d < M$, 
we claim 
$$D_\Gamma \geq {m\over  2}. \eqno{(5.6)}$$
To see this, suppose that there are $m/2$ good cubes. For each good cube $B_k(u)$, $\bar{B}_k(u)$ contains an edge $e\in R_\Gamma$ with $t(e)=1$,
so $e$ is a pivotal edge. Note that each $B_k(u)$ has at most $4^d$ neighbor cubes adjacent to $B_k(u)$,
so there are at least ${m\over 4^d2}$ pivotal edges. Therefore,  $K_\Gamma > {m\over 4^d2}$.
By (5.3) on $\{|R_\Gamma|\geq M K_\Gamma, |\hat{R}_\Gamma(k)|=m\}$,
$$|R_\Gamma|\geq M K_\Gamma \geq {Mm\over 4^d 2}\geq  {M|\hat{R}_\Gamma(k) |\over 4^d2}> 
{|\hat{R}_\Gamma (k)| k^d}.\eqno{(5.7)}$$
The contradiction of  (5.3) and (5.7) 
tells us that  (5.6) holds.\\

 By this observation and (5.4), we take $2(4k)^d  < M$ to obtain
$$P(|R_\Gamma|\geq MK_\Gamma)=\sum_{m\geq Mt/(2k^d)} P(|R_\Gamma|\geq Mt, |\hat{R}_\Gamma(k)|=m,D_\Gamma \geq m/2).\eqno{(5.8)}$$
Now we fix $\hat{R}_\Gamma$ to have
$$P(|R_\Gamma|\geq MK_\Gamma)=\sum_{m\geq Mt/(2k^d)}\sum_{\kappa_m} P(|R_\Gamma|\geq MK_\Gamma, \hat{R}_\Gamma(k)=\kappa_m ,D_\Gamma \geq m/2),\eqno{(5.9)}$$
where $\kappa_m$ is a fixed connected vertex set with $m$ vertices, and the second sum in (5.9) takes over all possible
such $\kappa_m$.
For each fixed $\hat{R}_\Gamma(k)=\kappa_m$, there are at most $m\choose i$ choices for these $i$,  $i=m/2,...,m$,
 bad cubes, so by 
(5.5)
$$P(|R_\Gamma|\geq MK_\Gamma, \hat{R}_\Gamma(k)=\kappa_m ,D_\Gamma \geq m/2)\leq C_1 m{m\choose m/2} \exp(-C_2 k m/2).\eqno{(5.10)}$$
Substitute the upper bound of (5.10) for each term of the sums in (5.9) to obtain
$$P(|R_\Gamma|\geq MK_\Gamma)\leq \sum_{m\geq Mt/(2k^d)}\sum_{\kappa_m}C_1m{m\choose m/2 } \exp(-C_2 k m/2).\eqno{(5.11)}$$
As we mentioned, $\hat{R}_\Gamma$ is connected, so there are at most $(4)^{dm}$ choices for  
all possible $\kappa_m $.
With this observation and (5.11) we have
$$P(|R_\Gamma|\geq MK_\Gamma)=\sum_{m\geq Mt/(2k^d)}(4)^{dm}m{m\choose m/2} \exp(-C_2 k m/2)\leq C_1 \sum_{m\geq Mt/k^d}
m[4^d 2 \exp(-C_2 k/2)]^m.\eqno{(5.12)}$$
We choose $k$ large to make
$$ 4^d 2 \exp(-C_2 k/2)<1/2.$$
 By (5.12), there are $C_3=C_3(p,d)$ and $C_4=C_4(p,d)$ such that
$$P(|R_\Gamma|\geq MK_\Gamma)\leq C_3 \exp(-C_4 t).\eqno{(5.13)}$$
Therefore, by (5.2) note that there are at most $t^{2d}$ vertices on $\Gamma$, so there exists $C_5=C_5(p,d)$
such that
$$E|R_\Gamma|=C_3 t^{2d} \exp(-C_4 t)+MCt\leq C_5 t.\eqno{(5.14)}$$
Theorem 2 follows from (5.14).

\newpage
\begin{center}
{\bf \large References}
\end{center}
Alexander, K (1993) A  note on some rated of convergence in first passage percolation. Ann. Appl. Probab. {\bf 3}
81-91.\\
Alexander, K. (1996) Approximation of subadditive functions and convergence rates in limiting-shape results.
Ann. Probab. {\bf 25} 30-55.\\
Grimmett, G. (1999) Percolation. Springer, Berlin.\\
Kardar, D. A. (1985) Roughening by impurities at finite temperatures. Phys. Rev.. Lett. {\bf 55} 2923-2923.\\
Kardar, D. A., Parisi, G. and Zhang, Y.C. (1986) Dynamic scaling of growing interfaces. Phys. Rev. Lett. {\bf 56} 889-892.\\
Kardar, D. A. and Zhang, Y.C. (1987) Scaling of directed polymers in random media. Phys. Rev. Lett. {\bf 56} 2087-2090.\\ 
 Hammersley J.M. and  Welsh D. J. A.(1965),
First-passage percolation, subadditive processes,
stochastic networks and generalized renewal theory,
in Bernoulli, Bayse, Laplace Anniversary Volume,
J. Neyman and L. LeCam eds., 61--110, Springer, Berlin.\\
Kesten, H. (1986), Aspects of first-passage percolation, Lecture Notes in
Mathematics {\bf 1180}, Springer, Berlin.\\
Kesten, H. (1993) On the speed of convergence in first passage percolation. Ann Appl. Probab. {\bf 3} 296-338.\\
Kesten, H and Zhang , Y (1990) The probability of a large finite clusters in supercritical Bernoulli percolation. Ann.
Probab. {\bf 18}, 537-555.\\
Kesten, H. and Zhang, Y. (1996) A  central limit theorem  for critical first passage percolation in two dimensions. PTRF {\bf 107} 137-160.\\
Kim, J. M. and Kosterlitz, M. (1989) Growth in a restricted solid on solid model. Phys. Rev. Letter {\bf 62} 2289-2292.\\
Krug, J. and Spohn, H. (1991) Kinetic roughening of growing surfaces. In solids Far from
equilibrium: Growth, Morphology Defects (C. Godreche, ed) 497-582. Cambridge Univ. Press.\\
Natterman, T. and Renz, W. (1988) Interface roughening due to random impurities at
low temperatures. Phys. Rev. B {\bf 38} 5184-5187.\\
Smythe R.T. and Wierman J. C.(1978),
First Passage Percolation on the Square Lattice,
Lecture Notes in Mathematics {\bf 671}, Springer, Berlin.\\
Wolf, D. and Kertesz, J. (1987) Surface with exponents for three and four dimensional Eden growth. Europhys Lett. {\bf 1} 651-656.\\

\noindent
Yu Zhang\\
Department of Mathematics\\
University of Colorado\\
Colorado Springs, CO 80933\\
email: yzhang@math.uccs.edu\\

\end{document}